\documentclass[10pt, a4paper, final]{amsart}
\pdfoutput=1
\usepackage{amstext, amscd, amsthm, bm, amssymb, dsfont}
\usepackage{amsmath}
\usepackage[english]{babel}
\usepackage{booktabs}
\usepackage{csquotes}
\usepackage[style=ieee-alphabetic,url=false,isbn=false,giveninits=true,backend=biber,maxnames=5,maxalphanames=5]{biblatex}
\bibliography{../bib}

\usepackage[utf8]{inputenc}

\usepackage{color}
\usepackage{datetime}
\usepackage{enumitem}
\usepackage{graphicx}
\usepackage{wasysym}
\usepackage{mathtools}
\usepackage[scr=rsfso]{mathalfa}
\usepackage{microtype}
\usepackage{textcomp}
\usepackage{hyperref}
\DeclareSymbolFont{extraup}{U}{zavm}{m}{n}
\DeclareMathSymbol{\varheart}{\mathalpha}{extraup}{86}

\usepackage{savesym}
\usepackage{stmaryrd}
\usepackage{textpos}
\savesymbol{iint}
\restoresymbol{TXF}{iint}
\usepackage{tikz}
\usetikzlibrary{babel}
\usepackage{tikz-cd}

\usepackage[backgroundcolor=blue!40, linecolor=lightgray, bordercolor=lightgray, textsize=tiny]{todonotes}

\usepackage[notref]{showkeys}

\newcommand{\R}             {\mathbb R} 
    
\newcommand{\N}             {\mathbb N}
\newcommand{\Z}             {\mathbb Z}
\newcommand{\Q}             {\mathbb Q}
\newcommand{\C}             {\mathbb C}

\newcommand{\cat}[1]{\pmb{\mathrm{#1}}}

\usepackage{xparse}
\renewcommand{\mathbf}[1]{\ensuremath{\pmb{\mathrm{#1}}}}

\DeclareDocumentCommand{\Ddeath}{ g }{
  \ensuremath{
    \mathrm{D}^{\dagger\IfNoValueF{#1}{,#1}}
  }
}
\DeclareDocumentCommand{\D}{  }{
  \ensuremath{
    \mathrm{D}
  }
}

\DeclareDocumentCommand{\Bdeath}{ g }{
  \ensuremath{
    \mathbf{B}_{\IfNoValueF{#1}{#1}}^{\dagger}
  }
}
\DeclareDocumentCommand{\Brigdeath}{ m g }{
  \ensuremath{
    \mathbf{B}_{\mathrm{rig}, #1}^{\dagger\IfNoValueF{#2}{,#2}}
  }
}

\DeclareDocumentCommand{\Drigdeath}{ g }{
  \ensuremath{
    \mathrm{D}_{\mathrm{rig}}^{\dagger\IfNoValueF{#1}{,#1}}
  }
}
\DeclareDocumentCommand{\Ars}{ g g }{
  \ensuremath{
    \mathbf{A}_{\IfNoValueF{#1}{#1}}^{\IfNoValueF{#2}{#2}}
  }
}

\DeclareDocumentCommand{\Brs}{ g g }{
  \ensuremath{
    \mathbf{B}_{\IfNoValueF{#1}{#1}}^{\IfNoValueF{#2}{#2}}
  }
}

\DeclareDocumentCommand{\Drs}{ g }{
  \ensuremath{
    \mathrm{D}^{\IfNoValueF{#1}{#1}}
  }
}

\usetikzlibrary{arrows,calc,matrix,decorations.pathmorphing}
\tikzset{sra/.code=
  {\pgfkeysalso{
      /tikz/decorate,
      /tikz/decoration={snake,amplitude=1pt,
        segment length=5pt, pre=lineto,
        pre length=4pt, post=lineto, post length=8pt}
    }}
}

\newcommand*{\ra}[1][]{\raisebox{-1pt}{\begin{tikzpicture}%
        \draw[xscale=1,thin,shorten >=3pt, >=stealth, ->,font=\scriptsize] (0,0)%
                 node{\hspace*{-2pt}} -- (0.5,0) node[above] {#1}%
                  -- node{} (1,0);\end{tikzpicture}}\penalty1000\relax}

\renewcommand*{\mapsto}{\raisebox{-1pt}{\begin{tikzpicture}%
       \draw[xscale=1,thin,shorten >=3pt, >=stealth, |->] (0,0)%
                node{\hspace*{0pt}}%
                -- node{} (1,0);\end{tikzpicture}}\penalty1000\relax}

\makeatletter
\newcommand{\osetwv}[2]{
  {\mathop{#2}\limits^{\vbox to -4\ex@{\kern-\tw@\ex@
   \hbox{\hspace{-0.2em}\scriptsize #1}\vss}}}}
\makeatother

\makeatletter
\newcommand{\oset}[2]{
  {\mathop{#2}\limits^{\vbox to -5\ex@{\kern-\tw@\ex@
   \hbox{\hspace{-0.25em}\scriptsize$#1$}\vss}}}}
\makeatother

\tikzset{commutative diagrams/.cd,
                     arrow style=tikz,
                     diagrams={>=stealth', line width=0.5pt}}

\DeclareRobustCommand{\gobblefive}[5]{}

\DeclarePairedDelimiter\abs{\lvert}{\rvert}%
\DeclarePairedDelimiter\norm{\lVert}{\rVert}%
\DeclarePairedDelimiter{\ceil}{\lceil}{\rceil}

\makeatletter
\let\oldabs\abs
\def\abs{\@ifstar{\oldabs}{\oldabs*}}
\let\oldnorm\norm
\def\norm{\@ifstar{\oldnorm}{\oldnorm*}}
\makeatother

\DeclareMathOperator{\Hom}{Hom}

\DeclareMathOperator{\Map}{Map}

\DeclareMathOperator{\image}{im}

\DeclareMathOperator{\id}{id}

\DeclareMathOperator{\sgn}{sgn}

\DeclareMathOperator{\coord}{coord}
\DeclareMathOperator{\CE}{\ensuremath{\mathfrak C}}
\DeclareMathOperator{\tw}{tw}
\DeclareMathOperator{\ad}{ad}
\DeclareMathOperator{\Ad}{Ad}
\DeclareMathOperator{\Tr}{Tr}

\usepackage{cleveref}

\crefname{none}{}{}
\crefname{(none)}{}{}
\creflabelformat{(none)}{(#2#1#3)}
\crefname{diagram}{diagram}{diagrams}
\creflabelformat{diagram}{(#2#1#3)}

\newtheoremstyle{mythm}
  {}
  {}
  {\itshape}
  {}
  {\bfseries\scshape}
  {.}
  { }
  {}
  
\newtheoremstyle{myprop}
  {}
  {}
  {\itshape}
  {}
  {\bfseries}
  {.}
  { }
  {}

\newtheoremstyle{myrmk}
  {}
  {}
  {}
  {}
  {\bfseries}
  {.}
  { }
  {}

\theoremstyle{mythm}
\newtheorem{theorem}{Theorem}[section]
\newtheorem*{theorem*}{Theorem}

\theoremstyle{myprop}
\newtheorem{lemma}[theorem]{Lemma}
\newtheorem{corollary}[theorem]{Corollary}

\newtheorem{proposition}[theorem]{Proposition}

\theoremstyle{myrmk}
\newtheorem{remark}[theorem]{Remark}
\newtheorem{example}[theorem]{Example}

\newtheorem{definition}[theorem]{Definition}
\newcommand{\rom}[1]{\uppercase\expandafter{\romannumeral #1\relax}}

\makeindex

\author{Oliver Thomas}
\address{Mathematisches Institut, Im Neuenheimer Feld 205, D-69120 Heidelberg}
\email{othomas@mathi.uni-heidelberg.de}
\title{Duality for $K$-analytic Cohomology}

\begin{document}
\renewcommand{\ref}[1]{\cref{#1}}
\begin{abstract}
  We prove a duality result for the analytic cohomology of Lie groups over
  non-archimedean fields acting on locally convex vector spaces.
\end{abstract}
\maketitle
\tableofcontents
\section*{Introduction}

Let $K$ be a non-archimedean complete field and $G$ a $K$-analytic group,
i.\,e., a group object in the category of $K$-analytic manifolds. Let
furthermore $V$ be an analytic representation of $V$ and $C^\bullet(G,V)$ the
complex of analytic inhomogeneous cochains of $G$ with coefficients in $V$. Its
cohomology is called the $K$-analytic cohomology of $G$ with coefficients in
$V$. If $K=\Q_p$ then Lazard showed that this is just continuous cohomology,
cf.~\cite[V.(2.3.10)]{MR0209286}. But if $K\neq\Q_p$ then this is no longer the
case and the $K$-analytic cohomology differs from the $\Q_p$-analytic
cohomology.

Assume for a moment that $V$ is of finite dimension over $K$. Then for $d=\dim
G$, we show the existence of a quasi-isomorphism \[ C^\bullet(G,V')
  \oset{\cong}{\ra} C^\bullet(G,V)'[-d],\] where $(-)'=\Hom_K(-,K)$. If $V$ is
not of finite dimension, then the functional analysis regrettably gets more
complicated: We still get a morphism \[ C^\bullet(G,V_b') \ra\Hom_K(
  C^\bullet(G,V),K)[-d]\] where $V_b'$ denotes the strong dual of $V$, but we
need additional requirements for this morphism to be a quasi-isomorphism. For
example, the Hahn-Banach theorem only holds for certain subclasses of
non-archimedean fields -- so taking the continuous dual is not always an exact
functor. The precise statement of our main
\ref{prop:duality-morphism-analytic-cohomology} includes more assumptions for
these kinds of reasons, which we will explain in the first few sections.

Strategically, the proof of the duality result is charmingly straight-forward:
Using \cite{MR3352825}, we compare analytic cohomology with Lie algebra
cohomology. Hazewinkel (cf.~\cite{MR0276285}) showed a duality result for Lie
algebra cohomology and plucking both results together then yields the result.

Technically, things are more complicated. The van Est comparison between
analytic cohomology and Lie algebra cohomology only yields an isomorphism of
cohomology groups when the underlying group is sufficiently connected --
something that of course isn't the case for non-archimedean ground fields.
Showing that the duality of the Lie algebra cohomology correctly identifies the
subspaces stemming from analytic cohomology is one issue, taking care of
multiple topological subtleties not present in the archimedean world another.

The main application we had in mind developing these results concerns Lubin-Tate
$(\varphi,\Gamma)$-modules as they appear in Iwasawa theory and especially the
Herr complex used to express its cohomology. Dualising this complex should (at
least morally) yield another Herr complex -- but we run into all sorts of
topological issues. Therefore, we can only construct a natural comparison
morphism. Whether it is a quasi-isomorphism is unknown to us.

\subsection*{Acknowledgements}

This article is based upon parts of its author's PhD thesis
(cf.~\cite{thomas:on-analytic-and-iwasawa-cohomology}). We thank Otmar Venjakob
for valuable discussions and insights.

\section{Some Functional Analysis}

We want to briefly recall some notions of non-archimedean functional analysis.
We refer the reader to \cite{MR0219078,MR1869547,MR3685952} for details. An
excellent overview can also be found in \cite{MR1664230}. In this section, we
fix a complete non-archimedean field $K$ with valuation ring $\mathcal O_K$.

\subsection{Foundations}

\begin{definition}
  $K$ is called spherically complete, if every decreasing sequence of closed balls has
  a non-empty intersection.
\end{definition}

\begin{example}
  Every locally compact field is spherically complete. $\C_p$, the completion of
  an algebraic closure of $\Q_p$, is not spherically complete.
\end{example}

\begin{definition}
  A lattice $L$ in a $K$-vector space $V$ is an $\mathcal O_K$-submodule of $V$,
  which satisfies \[ V = \bigcup_{\lambda\in K}\lambda L.\]
\end{definition}

\begin{definition}
  We call a topological $K$-vector space \emph{locally convex} (or an LCVS),
  if it has a neighbourhood basis of lattices. 
\end{definition}

\begin{remark}
  Note that a subset $M$ of a $K$-vector space is an $\mathcal O_K$-module if
  and only if for all $m,m'\in M$ and all $\lambda,\mu$ with
  $\abs{\lambda},\abs{\mu}\leq 1$ also $\lambda m+\mu m'\in M. $ This is the
  analogy to the usual notion of convexity. Requiring $\lambda m +
  (1-\lambda)m'\in M$ regrettably does not suffice.
\end{remark}

\begin{remark}
  Let $V$ be a $K$-vector space. For every lattice $L$ in $V$, there is an
  attached seminorm $p_L$ defined by \[ p_L(v) = \inf_{\lambda \in K,
      a\in\lambda L} \abs{\lambda}.\] Conversely, for a seminorm $p\colon V\ra
  \R$ and $\varepsilon> 0$ we can define a lattice \[ V_p(\varepsilon) = \{v\in
    V\mid p(v)<\varepsilon\}.\] These constructions are inverse to one another
  in the following sense: For a family of seminorms $(p_i)_i$, the coarsest
  topology on $V$ such that all $p_i$ are continuous is the locally convex
  topology generated by the lattices $(V_{p_i}(\varepsilon))_{i,\varepsilon}$.
  Conversely, if $V$ is locally convex, the topology on $V$ is the coarsest
  topology, such that all $(p_L)_{L}$ are continuous, where $L$ ranges over the
  open lattices in $V$. We refer to \cite[section I.4]{MR1869547} for details.
\end{remark}

\begin{definition}
  A subset $B$ of an LCVS $V$ is called bounded, if for any open lattice $L$ in
  $V$ there is a $\lambda\in K$ such that $B\subset\lambda L.$
\end{definition}

\begin{proposition}
  Every quasi-compact subset $C$ of an LCVS $V$ is bounded.
  \begin{proof}
    Let $L$ be an open lattice. By assumption, $V=\bigcup_{\lambda\in K} \lambda
    L$, so finitely many $\lambda_1 L,\dots, \lambda_n L$ cover $C$. We can
    assume that none of the $\lambda_i$ lie in $\mathcal O_K.$ Then $C\subseteq
    \lambda_1\cdots\lambda_n L.$
  \end{proof}
\end{proposition}

\begin{remark}
  If $K$ is not locally compact, an LCVS over $K$ does not have non-trivial
  compact $\mathcal O_K$-submodules.
\end{remark}

\begin{definition}
  Let $V$ be an LCVS. We call $V$ bornological, if a $K$-linear map $V\ra W$ of
  LCVS is continuous if and only if it respects bounded subsets. $V$ is called
  barrelled, if every closed lattice is open. 
\end{definition}

\subsection{Dual spaces}

\begin{definition}
  Let $V,W$ be LCVS. We denote the set of continuous $K$-linear maps from $V$ to
  $W$ by $\mathscr L(V,W).$ For bounded subsets $B\subseteq V$ and open subsets
  $U\subseteq W$ we denote by $L(B,U)\subseteq\mathscr L(V,W)$ those continuous
  linear maps which map $B$ into $U$. The families
  \begin{gather*}
    \left\{ L(S,U) \mid S\subseteq V\text{ a single point, } U\subseteq W\text{ open} \right\}\\
    \left\{ L(C,U) \mid C\subseteq V\text{ compact, } U\subseteq W\text{ open} \right\}\\
    \left\{ L(B,U) \mid B\subseteq V\text{ bounded, } U\subseteq W\text{ open} \right\}
  \end{gather*}
  generate locally convex topologies on the space $\mathcal{L}(V,W)$ of
  continuous linear maps from $V$ to $W$, which are called the \emph{weak},
  \emph{compact-open}, and \emph{strong} topology respectively. The
  corresponding LCVS will be denoted by $\mathscr L_s(V,W), \mathscr
  L_c(V,W),$ and $\mathscr L_b(V,W)$.
\end{definition}

\begin{remark}
  The weak topology is coarser than the compact-open topology, which in turn is
  coarser than the strong topology. 
\end{remark}

\begin{remark}
  \label{rmk:co-top-pathologies}
  Denote by $\cat{T}$ the category of Hausdorff topological spaces. (A
  variant of this remark also holds in the non-Hausdorff case.) For topological
  spaces $X,Y$ we denote by $[X,Y]$ the set $\Hom_{\cat{T}}(X,Y)$ endowed with
  the compact-open topology. It is an easy exercise to check that for
  topological spaces $X,Y,Z$ there is a well-defined map \[
    \Hom_{\cat{T}}(X\times Y, Z)\ra\Hom_{\cat{T}}(X,[Y,Z])\] sending $f$ to \[x
    \mapsto \left( y\mapsto f(x,y) \right).\] However, the obvious candidate for
  an inverse \[\Hom_{\cat{T}}(X,[Y,Z])\ra \Hom_{\cat{Set}}(X\times Y, Z),\]
  sending $f$ to \[ (x,y)\mapsto f(x)(y),\] in general does not yield continuous
  maps! Formally speaking, not every topological space is exponentiable. In our
  setting, we would have a bijection if $Y$ was locally compact, and locally
  compact spaces are the largest class for which this holds for all spaces $X$
  and $Z$. As LCVS are only locally compact if they are finite dimensional, we
  cannot use the adjointness properties of the compact-open topology. 
  In fact, there is mostly no reason to look at the compact-open topology at
  all. Considering linear maps, the strong topology plays the same role, but better.
\end{remark}

\begin{proposition}[Hahn-Banach]
  \label{prop:hahn-banach}
  If $K$ is spherically complete, $V$ a LCVS and $W$ a linear subspace of $V$
  endowed with the subspace topology. Then every continuous linear map $W\ra K$
  extends to a continuous linear map $V\ra K.$
  \begin{proof}
    \cite[proposition 9.2, corollary
    9.4]{MR1869547}.
  \end{proof}
\end{proposition}
 
There is also the following version of the Hahn-Banach theorem for LCVS of
countable type.

\begin{definition}
  \label{def:countable-type}
  An LCVS $V$ is said to be of countable type, if for every continuous seminorm
  $p$ on $V$ its completion $V_p$ at $p$ has a dense subspace of countable
  algebraic dimension.
\end{definition}

\begin{proposition}
  \label{prop:hahn-banach-countable-type}
  Let $V$ be an LCVS of countable type and $W$ a sub-vector space endowed with
  the subspace topology. Then every continuous linear map $W\ra K$
  extends to a continuous linear map $V\ra K.$
  \begin{proof}
   \cite[corollary 4.2.6]{MR2598517} 
  \end{proof}
\end{proposition}
\begin{definition}
  We say that Hahn-Banach holds for an LCVS $V$, if $K$ is spherically complete
  or $V$ is of countable type.
\end{definition}
\begin{remark}
  \label{rmk:countable-type-constructions}
  Spaces of countable type are stable under forming subspaces, linear images,
  projective limits, and countable inductive limits, cf.~\cite[theorem
  4.2.13]{MR2598517}.
\end{remark}

\begin{proposition}
  \label{prop:strong-top-adjointness}
  Let $f\colon V\times W\ra X$ be a (jointly) continuous bilinear map of LCVS.
  Then it induces a continuous map $f\colon V\ra\mathscr L_b(W,X)$.
  \begin{proof}
    As a jointly continuous map is also separately continuous, we have a
    well-defined map $f\colon V\ra \mathscr L(W,X).$ We only need to show that
    it is continuous with respect to the strong topology. For this, let
    $B\subseteq W$ be bounded and $M\subseteq X$ an open lattice. We need to
    show that the set $T$ of those $w\in W$ such that $f(w,B)\subseteq M$ is
    open. Let $w\in T$ and $b\in B$. By separate continuity we get open lattices
    $w\in L, b\in L'$ such that $f(L\times L')\subseteq M.$ As $M$ is bounded,
    there exists $\lambda\in K$ with $B\subseteq \lambda L'.$ Then \[
      f(w+\lambda^{-1}L, B) = f(w,B) + f(L,\lambda^{-1}B) \subseteq f(w,B) +
      f(L,L')\subseteq M.\]
  \end{proof}
\end{proposition}

\begin{proposition}[Banach-Steinhaus]
  \label{prop:banach-steinhaus}
  Let $V,W$ be LCVS. If $V$ is barrelled, then every bounded subset $H\subseteq \mathscr
  L_s(V,W)$ is equicontinuous, i.\,e., for every open lattice $L'\subseteq W$
  there exists an open lattice $L\subseteq V$ such that $f(L)\subset L'$ for
  every $f\in H.$
  \begin{proof}
   \cite[proposition 6.15]{MR1869547} 
  \end{proof}
\end{proposition}

\begin{proposition}
  Let $G$ be a locally compact topological group and $V$ a barrelled LCVS.
  Assume that $G$ acts via linear maps on $V$. Then \[G \times V\ra
    V\] is continuous if and only if it is separately continuous.
  \begin{proof}
    It is clear that a continuous group action is separately continuous.
    
    Let $U\subset V$ be an open lattice and $g\in G, v\in V, gv\in U.$ Let $H$ be a compact
    neighbourhood of $g$ with $Hv\in U,$ which exists by local compactness of
    $G$ and separate continuity of the group action.

    Consider the set $M=\{ h\cdot -\mid h\in H\}$ of continuous linear maps
    $V\ra V$. We want to show that it is bounded in the topology of pointwise
    convergence on $\Hom_{\mathrm{cts}}(V,V).$ For this matter, take $w\in V$,
    $S\subset V$ an open lattice, and denote by $L$ those continuous linear maps
    $V\ra V$ which map $w$ into $S$. We need to show that there exists
    $\lambda\in K$ with $M\subseteq \lambda L.$ As $H$ is compact, so is
    $Hw\subseteq V,$ hence there exists $\lambda\in K$ such that
    $Hw\in\lambda^{-1} S$, i.\,e., $M\subseteq \lambda L.$

    \Cref{prop:banach-steinhaus} now shows the existence
    of an open lattice $L'$ such that $HL'\subseteq U,$ or in other words,
    $H\times L'\subseteq \mathrm{mult}^{-1}(U).$
  \end{proof}
\end{proposition}

\begin{definition}
  The dual space of an LCVS $V$ is the vector space of continuous $K$-linear
  functions $V\ra K$ and will be denoted by $V'.$

  We denote by $V'_s$ the dual space equipped with the weak topology, which is
  the topology of pointwise convergence.

  $V'_c$ will denote the dual space equipped with the compact-open topology,
  which is also the topology of uniform convergence on compact subsets.

  The strong dual will be denoted by $V'_b$ and is defined as the topology of
  uniform convergence on bounded subsets of $V$. 
\end{definition}
\begin{remark}
  Note that for both the weak and strong duals, the dual of a direct sum of LCVS
  is the product of its duals. However, only for the strong dual is the dual of
  a product of LCVS the sum of its duals.

  Note that by \cite[lemma 6.4]{MR1869547}, $V'_c$ can be defined as the
  coarsest topology on $V'$ such that for every quasi-compact $K\subset V$ the
  map
  \begin{gather*}
    V' \ra \R\\
    v' \mapsto \sup_{v\in C} \abs{v'(v)}_K
  \end{gather*}
  is continuous.

\end{remark}

\subsection{Analyticity}

\begin{definition}
  \label{def:analytic-map}
  Let $E$ be a normed $K$-vector space and $V$ a LCVS. A formal sum \[ f=
    \sum_{n\in\N_0} f_n\] of continuous functions $f_n\colon E\ra V$ which are
  homogeneous of degree $n$ (i.\,e., $f_n(\lambda x)=\lambda^n f_n(x)$ for all
  $n\in\N_0,\lambda\in K, x\in E$) is called a convergent power series, if there
  exists an $R>0$ such that for every continuous seminorm $p\colon V\ra K$ the
  following holds: \[ \norm{(f_n)_n}_{p,R} = \sup_{n\in\N_0} \sup_{x\in E,
      \norm{x}\leq 1} R^n p(f_n(x)) < \infty.\]
  The supremum over all $R$ such that for every continuous seminorm $p$ we have
  $\norm{(f_n)_n}_{p,R}<\infty$ is called the radius of convergence of $f.$

  A map $\widetilde{f}\colon E\ra V$ is called analytic in $x\in E,$ if there
  exists a convergent power series $f_x$ such that for all $h\in E$ close enough
  to zero, we have an equality \[ \widetilde{f}(x+h) = f_x(h).\] It is called
  analytic, if it is analytic at every point.
  
  Let $M$ be an analytic Banach manifold over $K$ (i.\,e., $M$ is locally
  isomorphic to $K$-Banach spaces with analytic transition maps). A map
  $\widetilde{f}\colon M\ra V$ with values in $V$ is called locally analytic, if
  it is analytic in charts. The radius of convergence of $\widetilde{f}$ at $x$
  is the radius of the power series development at a local chart. It might be
  larger than the chart itself.
\end{definition}

\begin{lemma}
  \label{prop:radius-of-convergence-lower-semicts}
  Let $f\colon E\ra V$ be an analytic map from a normed vector space to a
  Hausdorff LCVS $V$. The map $r_f\colon E\ra\R_{>0}\cup\{\infty\}$, mapping a
  point to the radius of convergence of the power series development of $f$ at
  that point, is lower
  semi-continuous, i.\,e., for every $x\in E$ we have \[\liminf_{x'\to x}
    r_f(x')\geq r_f(x).\] Consequently, if $C\subseteq E$ is compact, then \[
    \inf_{x\in C} r_f(x) > 0.\]
  \begin{proof}
    A power series is analytic within its ball of convergence. It follows that
    the radius of convergence can at most increase. As lower semi-continuous maps
    attain their infimum in compact sets, the claim follows.
  \end{proof}
\end{lemma}

\begin{proposition}
  \label{prop:power-series-unique}
  The development as a power series is unique, i.\,e., if $\widetilde{f}\colon
  E\ra V$ is an analytic map between a normed $K$-vector space $E$ and a
  Hausdorff LCVS $V$ and if \[ \widetilde{f}(x+h)=\sum_{n\in\N_0}
    f^{(1)}_{x,n}(h) = \sum_{n\in\N_0} f^{(2)}_{x,n}(h) \] for sufficiently
  small $h$ with $f^{(i)}_{x,n}$ continuous and homogeneous of degree $n$, then
  $f^{(1)}_{x,n}=f^{(2)}_{x,n}$ for all $n$.
    \begin{proof}
      It suffices to show that if $\sum_n f_n$ is the zero function with $f_n$
      continuous and homogeneous of degree $n$ and $(f_n)_n$ convergent close to
      zero, then all $f_n=0$. Assume that $f_k\neq 0$. We can assume that $k$ is
      minimal with this property. Let $p$ be a continuous seminorm on $V$ with
      $p(f_k(x))> 0$. By replacing $x$ with $\lambda x$ for some $\lambda$ close
      to zero, we can assume that $p(f_k(x)) > p(f_{k+n}(x))$ for all $n>0$: By
      convergence of the power series, $\{p(f_{k+n}(x))\mid n\}\subseteq \R$ is
      bounded from above by some $R\in\R$. Choose now $\lambda\in K$ with
      $\abs{\lambda}<\max\{1, p(f_k(x))/R\}$, then it is easy to see that indeed
      $p(f_k(\lambda x))\neq 0$, $k$ is minimal with this property, and
      $p(f_k(\lambda x)) > p(f_{k+n}(\lambda x))$ for all $n>0.$

      But then $p(\sum_n f_n(x))=p(f_k(x))>0$, so $\sum_n f_n$ is not the zero
      function.
    \end{proof}
\end{proposition}

\subsection{Strictness}

\begin{definition}
  A linear map $V\ra W$ of LCVS is called \emph{strict}, if the induced map \[
    V/\ker f\ra\image f\] with the quotient topology on $V/\ker f$ and the subspace
  topology on $\image f$ is an isomorphism. 
\end{definition}

\begin{remark}
  Open linear maps are clearly strict, but strictness is remarkably bad behaved
  in general: Neither the sum nor the composition of strict maps needs to be
  strict again.
\end{remark}

\begin{definition}
  A sequence of LCVS \[0\ra A\ra B\ra C\ra 0\] is called exact if it is exact as
  a sequence of vector spaces and if the involved maps are all strict. 
\end{definition}
\begin{proposition}
  \label{prop:spherically-complete-implies-exact-dual}
  \label{prop:countable-type-exact-dual}
  Let \[0\ra A\ra B\ra C\ra 0\] be an exact sequence of LCVS. If Hahn-Banach
  holds for $B$, then the induced sequence of abelian groups \[0\ra C'\ra
    B'\ra A'\ra 0\] is also exact.
  \begin{proof}
    Let \[0\ra A\oset{\iota}{\ra} B\oset{\pi}{\ra} C\ra 0\] be exact. It is
    clear that \[0\ra C'\oset{\pi^*}{\ra} B'\] is exact. Let $f\colon B\ra K$ be
    in the kernel of $\iota^*$, i.\,e., $\iota A\subseteq\ker f.$ This induces a
    map $B/\iota A\ra K,$ which by strictness is a map $C\ra K.$

    It remains to show surjectivity of $\iota^*$, i.\,e., the existence of a
    map $\widetilde{f}$ such that the following diagram commutes: \[
      \begin{tikzcd}[row sep=0.1em]
        \iota A \arrow{r}{\cong} & A\arrow[r, "f"] & K\\
        \rotatebox[origin=c]{-90}{$\subseteq$} & & \\
        B \arrow[uurr, dashed, "\widetilde{f}"']& &\\
      \end{tikzcd}
    \]
    But this extension exists by \ref{prop:hahn-banach} or
    \ref{prop:hahn-banach-countable-type}.
  \end{proof}
\end{proposition}
\begin{lemma}
  \label{prop:matrices-are-strict}
  Let $V$ be an LCVS and $A\colon K^n\ra K^m$ a linear map. Then the induced
  map \[A\otimes_K V\colon V^{\oplus n}\ra V^{\oplus m}\] is strict.
  \begin{proof}
    Note that finite direct sums of LCVS coincide with their product. It is then
    clear that every component map $(A\otimes_K V)_i\colon V^{\oplus n}\ra V$ is
    open, so $A\otimes_k V$ is open as well.
  \end{proof}
\end{lemma}

\begin{definition}
  \label{def:frechet-space}
  A Fréchet space is an LCVS which is isomorphic to the projective limit of
  Banach spaces.
\end{definition}
\begin{remark}
  A space is Fréchet if and only if it is a complete LCVS whose topology is induced
  by a translation-invariant metric if and only if it is a Hausdorff topological
  $K$-vector space whose topology is induced by a countable family of semi-norms
  for which every Cauchy sequence converges. 
\end{remark}

\begin{proposition}[Open-mapping theorem]
  \label{prop:open-mapping-thm}
  Let $f\colon V\ra W$ be a continuous surjective linear map from a Fréchet
  space to a barrelled Hausdorff LCVS. Then $f$ is open.
  \begin{proof}
   \cite[proposition 8.6]{MR1869547} 
  \end{proof}
\end{proposition}

\begin{definition}
  \label{def:lf-space}
  An LCVS is called an LF-space, if it is the direct limit of a countable family
  of Fréchet spaces, the limit being formed in the category of locally convex
  vector spaces.
\end{definition}

\begin{remark}
  LF-spaces are Hausdorff. 
\end{remark}

\begin{proposition}[Open-mapping theorem for LF-spaces]
  \label{prop:open-mapping-thm-lf}
  Every continuous surjective linear map between LF-spaces is open. 
  \begin{proof}
    \cite[proposition 8.8]{MR1869547}  
  \end{proof}
\end{proposition}

\begin{proposition}
  If a continuous linear map $f\colon V\ra W$ between LF spaces has
  finite-dimensional cokernel, it is strict.
  \begin{proof}
    Note that this does not follow immediately from \ref{prop:open-mapping-thm-lf},
    as we do not know that $\image f$ is again LF.  

    Take finitely many independent vectors whose projection to
    the cokernel form a basis of the cokernel. Their span in $W$ will be called
    $X$. As $X$ is finite dimensional, it is especially also LF and hence
    so is $V\oplus X.$ The map
    \[
      \begin{tikzcd}
        V\oplus X\arrow{rr}{f\oplus\id}&& W
        \end{tikzcd}
    \] is then
    bijective, linear and continuous; thus it is an isomorphism by
    \ref{prop:open-mapping-thm} and $f$ hence an isomorphism onto $f(V)$.
  \end{proof}
\end{proposition}
\section{Analytic Actions of Lie Groups}

We continue with a fixed non-archimedean field $K$.

\begin{definition}
  A group object in the category of (finite-dimensional analytic)
  $K$-ma\-ni\-folds is called a Lie group over $K$.
\end{definition}

\begin{definition}
  \label{def:equianalytic}
  Let $G$ be a Lie group over $K$ and $V$ a separated LCVS. A continuous action
  $G\times V\ra V$ by continuous linear maps is called analytic, if every orbit
  map $g\mapsto gv$ is analytic. It is called equi-analytic, if it is analytic
  and the contragradient action on the dual space $G\times V'\ra V'$ is analytic
  with respect to the strong topology on $V'.$
\end{definition}

\begin{proposition}
  \label{prop:cts-eval-cts-contragradient}
  If a Lie group $G$ acts continuously on an LCVS $V$ and if the evaluation map
  $V'_b\times V\ra K$ is continuous, then the contragradient action $G\times
  V'_b\ra V'_b$ is also continuous.
  \begin{proof}
    Consider the following maps:
    \[
      \begin{tikzcd}
        G\times V\times V'_b\arrow{rr}{((-)^{-1},\id,\id)} && G \times V \times
        V'_b \arrow{rr}{(\mathrm{mult},\id)} && V \times V'_b \arrow{r}& K
      \end{tikzcd}
    \]
    The last map is just the evaluation function. The composite is now clearly
    continuous and by \ref{prop:strong-top-adjointness}, so is the induced map \[
      G\times V'_b\ra V'_b,\] which is the contragradient action.
  \end{proof}
\end{proposition}

We will spell out the following proposition in more detail than necessary to
show where the name \emph{equi-analytic} stems from.

\begin{proposition}
  \label{prop:fin-dim-is-equianalytic}
  An analytic action $G\times V\ra V$ is equi-analytic, if $V$ is of finite
  dimension.
  \begin{proof}
    By \ref{prop:cts-eval-cts-contragradient} we only need to show that every
    orbit map \[g \mapsto v'(g^{-1}\cdot -)\] is analytic. Fix $v'\in V'.$

    Considering a chart $\coord\colon U\ra K^d$ of a neighbourhood $U$ of $g$
    and $h$ close to the neutral element, \[ (gh)^{-1}v = \sum_{n\in\N_0}
      F_{g,v,n}(\mathrm{coord}(h))\] with $F_{g,v,n}\colon K^d\ra V$ continuous
    and homogeneous of degree $n$. Define \[F_{g,v',n}'(x)(v) =
      v'(F_{g,v,n}(x)).\] It suffices to show that \[ F'_{g,v',n}\colon U\ra V'_c\]
    is well-defined, continuous, homogeneous of degree $n$, and gives rise to a
    convergent power series, as then
    \begin{align*}
      \sum_n F'_{g,v',n}(\mathrm{coord}(h))(v) &= \sum_n v'(F_{g,v,n}(\mathrm{coord}(h)))\\
                                                          &= v'(\sum_n F_{g,v,n}(\mathrm{coord}(h)))\\
                                                          &= v'((gh)^{-1}v)\\
                                                          &= ((gh) v')(v).
    \end{align*}
    Note first that by linearity of $v',$ indeed $F'_{g,v',n}$ is homogeneous of
    degree $n$. Using \ref{prop:power-series-unique}, we see that
    $F'_{g,v',n}(h)$ is $K$-linear. The same argument that resulted in
    \ref{prop:cts-eval-cts-contragradient} also shows that $F'_{g,v',n}$ is
    continuous. It remains to show that $(F'_{g,v',n})_n$ is convergent with
    respect to the strong topology, i.\,e., we need to show that there
    exists an $R>0$ such that for every bounded set $B\subset V$ we have that \[
      \sup_{n\in\N_0} \sup_{x\in K^d,\norm{x}\leq 1} \sup_{v\in B} R^n
      \abs{F'_{g,v',n}(x)(v)}<\infty,\] which by definition of $F'_{g,v',n}$ is
    equivalent to \[ \tag{$\star$}\sup_{n\in\N_0} \sup_{x\in K^d,\norm{x}\leq 1} \sup_{v\in
        B} R^n \abs{v'(F_{g,v,n}(x))}<\infty.\] Analyticity of the group action
    on the other hand yields that for fixed $g\in G,v\in V$ we have an
    $R_{g,v}>0$ such that \[\sup_{n\in\N_0} \sup_{x\in K^d,\norm{x}\leq 1}
      \sup_{v''\in V'} R_{g,v}^n \abs{v''(F_{g,v,n}(x))}<\infty.\] If
    $B\subseteq V$ is compact, \ref{prop:radius-of-convergence-lower-semicts}
    yields the existence of $R_{g,B}>0$ such that
    \[\tag{$\blacklozenge$}\sup_{n\in\N_0} \sup_{x\in K^d,\norm{x}\leq
        1}\sup_{v\in B}
      \sup_{v''\in V'} R_{g,C}^n \abs{v''(F_{g,v,n}(x))}<\infty.\]
    We cannot directly deduce $(\star)$ from $(\blacklozenge)$, as we have no
    means of controlling the radius of convergence across different compact (or bounded)
    subsets. This homogeneity is what \emph{equi-analytic} alludes to.

    By \ref{prop:power-series-unique}, for any $u,w\in V$ and $\lambda\in K$ \[
      F_{g,u+\lambda w, n}(x) = F_{g,u,n}(x) + \lambda F_{g,w,n}(x).\] It
    follows that if $u\in V$ is in the linear subspace generated by $w_1,\dots,
    w_k\in V,$ then for the radii of convergence of the orbit maps $-\cdot u$
    and $-\cdot w_i$ we have the following estimate: \[ r_{-\cdot u}(g)\geq
      \min_i r_{-\cdot w_i}(g).\] If $V$ is generated by $v_1,\dots,v_n$ and
    $R=\min_i r_{-\cdot v_i}(g)$, then $R>0$ and \[ \sup_{n\in\N_0} \sup_{x\in
        K^d,\norm{x}\leq 1} \sup_{v\in V}\sup_{v'\in V} R^n \abs{
        v'(F_{g,v,n}(x))} < \infty,\] which is more than enough to show
    $(\star).$
  \end{proof}
\end{proposition}
\begin{lemma}
  \label{prop:dual-is-functor}
  Let $\varphi$ be a continuous endomorphism of $V$. Then it induces a
  continuous map $\varphi'\colon V'_b\ra V'_b.$
  \begin{proof}
    We need to show that if $B\subseteq V$ is bounded and $U\subseteq K$ is
    open, then also $(\varphi')^{-1}(L(B,U))=L(\varphi(B),U)$ is open. But a
    continuous map clearly maps bounded sets to bounded sets. 
  \end{proof}
\end{lemma}

\begin{proposition}
  \label{prop:cts-after-analytic-is-analytic}
  Let $M$ be a Banach manifold and $V,W$ separated LCVS. If $f\colon M\ra
  V$ is analytic and $\varphi\colon V\ra W$ continuous and linear, then
  $\varphi\circ f\colon M\ra W$ is also analytic.
  \begin{proof}
    We can assume that $M$ is a Banach space. Let $x\in M$ be arbitrary and let
    $f_{x,n}\colon M\ra V$ be continuous maps, homogeneous of degree $n$, such
    the family $(f_{x,n})_n$ is a convergent power series and that for all $h$
    sufficiently close to zero we have an equality \[ f(x+h) = \sum_n
      f_{x,n}(h).\] It suffices to show that the family $(\varphi\circ
    f_{x,n})_n$ is a convergent power series. By continuity of $\varphi$, for
    every continuous seminorm $p$ on $W$ we can find a $\lambda_p\in\R$ such
    that \[ \rho(\varphi(y))\leq\lambda_p\norm{y}. \] Let $R$ be the radius of
    convergence of $(f_{x,n})_n$ and $p$ a continuous seminorm on $W$. Then \[
      \sup_{n\in\N_0} \sup_{h\in E,\norm{h}\leq 1} R^n p(\varphi(f_{x,n}(h))) \leq
      \lambda_p \sup_{n\in\N_0} \sup_{h\in E,\norm{h}\leq 1}
      R^n\norm{f_{x,n}(h)} < \infty\]
  \end{proof}
\end{proposition}
\section{Duality for Lie Algebras}

For the general theory of Lie algebras and Lie groups we refer to
\cite{MR1176100,MR979493}. In this section, we fix a complete non-archimedean
field $K$ of characteristic zero and a Lie group $G$ over $K$. We also consider
its attached Lie algebra $\mathfrak g$ with Lie bracket $[-,-].$ The adjoint
action of $G$ on $\mathfrak g$ by differentiating conjugation maps will be
denoted by $\Ad(-)$, the adjoint action of $\mathfrak g$ on itself given by
$x\mapsto [x,-]$ will be denoted by $\ad(-)$.

\begin{definition}
  For a $\mathfrak g$-module $M$ we define the Chevalley-Eilenberg
  complex \[\textstyle \CE^\bullet(\mathfrak g, V) = \Hom(\bigwedge^\bullet
    \mathfrak g, V)\] concentrated in non-negative degrees by considering the
  differential \[ d\colon \CE^n(\mathfrak g, V)\ra \CE^{n+1}(\mathfrak g, V)\]
  given by
  \begin{align*} df(x_1\wedge\dots\wedge x_{n+1}) &=
                                                    \sum_i (-1)^{i+1} x_i f(x_1\wedge\dots\wedge\widehat{x_i}\wedge\dots\wedge x+{n+1})\\
                                                  & \quad + \sum_{i<j}
                                                    (-1)^{i+j} f([x_i,x_j]\wedge
                                                    x_1\wedge\dots\wedge\widehat{x_i}\wedge\dots\wedge\widehat{x_j}\wedge\dots\wedge
                                                    x_{n+1}).
  \end{align*}
  As usual, $\widehat{x_i}$ means omitting $x_i$ etc.
\end{definition}

\begin{definition}
  Let $V$ be a $\mathfrak g$-module. Define $V^{\tw}$ as the vector space $V$
  with a $\mathfrak g$-action given by \[ x\cdot^{\tw}v=xv-\Tr(\ad(x))v,\] where
  $\Tr$ is the trace map.
\end{definition}

\begin{proposition}
  For $V\neq 0,$ $V^{\tw}=V$ if and only if $H^{\dim \mathfrak g}(\mathfrak g,
  K)\neq 0$. If $\mathfrak g$ is abelian or nilpotent, $V^{\tw}=V.$
  \begin{proof}
    {\cite[corollary 2]{MR0276285}}
  \end{proof}
\end{proposition}

In applications, this is very often the case.

\begin{proposition}
  \label{prop:compact-group-trivial-lie-twist}
  If $G$ is compact, then $V^{\tw}=V.$
  \begin{proof}
    As for a compact group $G$, the left and right Haar measures coincide,
    \cite[section III.3.16]{MR979493} implies that \[ \det\Ad g = 1 \] for all $g\in G.$
    By \cite[section III.4.5]{MR979493}, we see that for all $x$ in a neighbourhood of zero of
    $\mathfrak g$
    \[ \Ad(\phi(x))=\exp(\ad x),\]
    where $\phi$ is a local exponential map from this neighbourhood into $G$.
    Here, $\exp$ is the usual exponential map of $K$ extended to matrices.
    Applying the determinant, we see that \[1 = \det\exp(\ad x) = \exp(\Tr\ad
      x)\] so $\Tr\ad x=0$ in a neighbourhood of the identity. Choosing a basis
    of $\mathfrak g$ in this neighbourhood, we see that indeed \[ \Tr(\ad x)
      =0\] for all $x\in G$ and hence $V^{\tw}=V.$
  \end{proof}
\end{proposition}

\begin{remark}
  The argument of \ref{prop:compact-group-trivial-lie-twist} shows that if
  $\Tr\ad x=0$ for all $x\in\mathfrak g$, then $\det\Ad(g)=1$ for all $g$ in a
  neighbourhood of the identity. If $G$ is a connected Lie group over $\R$ or
  $\C$, then $\det\Ad g=1$ for all $g\in G$. \cite[section III.3.16]{MR979493}
  then implies that the left and right Haar measures of $G$ coincide.
\end{remark}

\begin{definition}
  For any natural number $n$ we denote by $\ceil{n}$ the ordered set
  $\ceil{n}=\{1,\dots,n\}.$ For an injective morphism of ordered sets
  $\phi\colon\ceil{k}\ra\ceil{d}$ there exists a unique morphism of ordered
  sets \[\phi^*\colon \ceil{d-k}\ra\ceil{d}\] such that
  $\ceil{d}=\image\phi\cup\image\phi^*.$
  We then define \[\sgn\phi=(-1)^{\sum_{i=1}^{d-k}\phi^*(i)}.\]
\end{definition}

\begin{remark}
  \label{rmk:shuffle-sign-identity}
  It is easy to see that for an injective morphism of ordered sets
  $\phi\colon\ceil{k}\ra\ceil{d}$ the following holds: \[
    \sgn(\phi)\cdot\sgn(\phi^*)=(-1)^{k(d-k)},\] cf.~\cite[lemma
  10.17]{thomas:cohomology-of-topologised-monoids}.
\end{remark}

\begin{proposition}
  \label{prop:hodge-star-adjoint}
  Let $M$ be a finite dimensional vector space with basis $e_1,\dots,e_d$. For
  an injective morphism of ordered sets $\phi\colon\ceil{k}\ra\ceil{d}$
  define \[ e_\phi = e_{\phi(1)}\wedge\dots\wedge
    e_{\phi(k)}\in\textstyle\bigwedge^k M.\] Also define the $K$-linear
  isomorphism \[ \star\colon\textstyle\bigwedge^k M\ra\bigwedge^{d-k} M\]
  given by \[\star e_\phi = \sgn(\phi^*) e_{\phi^*}.\] Then for any invertible
  endomorphism $A$ of $M$ the following holds: \[
    \det A\cdot ((A^{-1})^t\circ\star) = \star\circ A.\]
  \begin{proof}
    This is a straight-forward piece of linear algebra, but we could not find a
    reference for $k\neq 1$.

    Let $\phi,\psi\colon \ceil{k}\ra\ceil{d}$ be injective maps of ordered sets.
    For a matrix $A$ denote by $A_{\phi,\psi}$ the matrix with entries
    $(a_{\phi(i),\psi(j)})_{i,j\in\ceil{k}}$. Now a straight forward calculation
    (or \cite[proposition 9 in III.8.5]{MR1727844}) shows that for fixed $\psi$,
    we have \[ Ae_\psi = \sum_\phi (\det
      A_{\phi,\psi})e_\phi,\tag{$\spadesuit$}\] where $\phi$ ranges
    over the injective maps of ordered sets $\ceil{k}\ra\ceil{d}.$ We hence also
    get \[ (A^{-1})^t\star e_\psi = \sgn(\psi^*)\sum_{\phi^*} \det((A^{-1})^t_{\phi^*,\psi^*})
      e_{\phi^*},\] where $\phi^*$ ranges over the injective maps of ordered sets
    $\ceil{d-k}\ra\ceil{d}$. Applying $\star$ to ($\spadesuit$), we are reduced
    to showing \[ \sgn(\psi^*)\cdot\det(A)\cdot \det (A^{-1})^t_{\phi^*,\psi^*} = \sgn(\phi^*)\cdot\det
      A_{\phi,\psi}.\]
    For $k=1,$ this is precisely Cramer's rule.
    
    Generally, for a matrix $B$, the submatrix $B_{\phi,\psi}$ can be considered
    as a linear map from the span of $e_{\phi(1)},\dots,e_{\phi(k)}$ to the span
    of $e_{\psi(1)},\dots,e_{\psi(k)}.$ Denote this linear map by
    $B^{\mathrm{res}}_{\phi,\psi}.$
    Define $B^{\mathrm{ext}}_{\phi,\psi}$ via
    \begin{gather*}
      B^{\mathrm{ext}}_{\phi,\psi} e_{\phi(i)} = B^{\mathrm{res}}_{\phi,\psi}e_{\phi(i)},\\
      B^{\mathrm{ext}}_{\phi,\psi} e_{\phi^*(i)} = e_{\psi^*(i)}.
    \end{gather*}
    It is clear that $\det B^{\mathrm{ext}}_{\phi,\psi}=\varepsilon \det
    B_{\phi,\psi},$ with $\varepsilon=(-1)^{\sum_i \phi^*(i)+\psi^*(i)}$
    so
    \begin{align*}
      \varepsilon\det (B_{\phi,\psi})\cdot e_1\wedge\dots\wedge e_n &=
      B^{\mathrm{ext}}_{\phi,\psi} e_1\wedge\dots\wedge
      B^{\mathrm{ext}}_{\phi,\psi} e_n\\ &= \sgn(\phi^*)\cdot
      B^{\mathrm{res}}_{\phi,\psi}e_{\phi}\wedge e_{\psi^*}
    \end{align*}
    By anticommutativity of the exterior algebra, we see that indeed 
    \[B^{\mathrm{res}}_{\phi,\psi}e_{\phi}\wedge e_{\psi^*} = Be_{\phi}\wedge
      e_{\psi^*},\]
    so
    \[\varepsilon\cdot \det (B_{\phi,\psi})\cdot e_1\wedge\dots\wedge e_n=\sgn(\phi^*)\cdot
      Be_{\phi}\wedge e_{\psi^*}.\tag{$\clubsuit$}
    \]
    If $B$ is invertible, we can apply $B^{-1}$ and get
    \begin{align*}
      \varepsilon\det
      (B_{\phi,\psi})\cdot\det(B^{-1})\cdot e_1\wedge\dots\wedge e_n &=
                                                                  \sgn(\phi^*)\cdot e_{\phi}\wedge B^{-1}e_{\psi^*}\\
      & = \sgn(\phi^*)(-1)^{k(d-k)}
      B^{-1}e_{\psi^*}\wedge e_\phi.
    \end{align*}
    
    Applying ($\clubsuit$) to $(B^{-1})_{\psi^*,\phi^*}$, we find that
    \[ B^{-1}e_{\psi^*}\wedge e_\phi=\varepsilon'\cdot \sgn(\psi)\cdot
      \det((B^{-1})_{\psi^*,\phi^*})\cdot e_1\wedge\dots\wedge e_n \] with
    $\varepsilon'=(-1)^{\sum_i \phi(i)+\psi(i)}$. As clearly
    $\varepsilon\varepsilon' =1,$ we get \[ \det(B_{\phi,\psi})\cdot
      \det(B^{-1}) = \sgn(\phi^*)\cdot
      (-1)^{k(d-k)}\cdot\sgn(\psi)\cdot\det((B^{-1})_{\psi^*,\phi^*}),\] and
    using \ref{rmk:shuffle-sign-identity}, this becomes \[
      \sgn(\phi^*)\cdot\det(B_{\phi,\psi})\cdot \det(B^{-1}) = \sgn(\psi^*)\cdot
      \det((B^{-1})_{\psi^*,\phi^*}),\] which is exactly what we needed to show.
  \end{proof}
\end{proposition}

\begin{theorem}
  \label{prop:chevalley-eilenberg-duality}
  Let $V$ be a LCVS with a continuous action by $\mathfrak g$. If $\mathfrak g$
  if of dimension $d$, and if $\Tr(\ad(x))=0$ for all $x\in\mathfrak g$, then
  there is a $G$-equivariant functorial isomorphism of complexes \[
    \CE^\bullet(\mathfrak g, V') \cong \CE^\bullet(\mathfrak g, V)'[-d].\]
  \begin{proof}
    In \cite{MR0276285}, Hazewinkel shows (without the restriction
    $\Tr\circ\ad=0$) that as abstract vector spaces \[ \CE^\bullet(\mathfrak g,
      (V^{\tw})^*) \cong \CE^\bullet(\mathfrak g, V)^*[-d],\] where
    $(-)^*=\Hom_K(-,K).$ While it is easy to check that the isomorphism respects
    continuous maps, it is not immediate at all that it is $G$-equivariant. The
    proof itself is a brutal calculation.

    Choose a basis $(e_i)_i$ of $\mathfrak g$ and define the star operator as in
    \ref{prop:hodge-star-adjoint}. Hazewinkel's isomorphism stems from the
    following pairing:
    \begin{gather*}
      \textstyle \langle -,- \rangle\colon \Hom_K(\bigwedge^k\mathfrak g,V')
      \times \Hom_K(\bigwedge^{d-k} \mathfrak
      g, V) \ra K\\
      (a,b)\mapsto\langle a,b \rangle=\sum_{\phi} a(e_{\phi})(b(\star e_{\phi}))
    \end{gather*}
    We need to show that \[ \langle ga,b \rangle = \langle a,g^{-1}b \rangle\]
    for all $g\in G.$ Write $A$ for $\Ad(g)$. Then \[ \langle gx,y \rangle =
      \sum_{\phi} (g.x)(e_{\phi})(y(\star e_{\phi})) = \sum_{\phi}
      (g(x(A^{-1}e_{\phi})))(y(\star e_{\phi}))=\sum_\phi
      x(A^{-1}e_{\phi})(g^{-1}y(\star e_{\phi})) \] and \[ \langle x,g^{-1}y
      \rangle = \sum_\phi x(e_\phi)(g^{-1} y(A\star e_\phi)). \]
    In both cases, $\phi$ runs over the injective increasing maps
    $\ceil{k}\ra\ceil{d}$.
    
    By considering the finite dimensional subspace of $V$ generated by all
    $y(\star e_\phi)$ and their images under $g^{-1}$, we can consider a finite
    dimensional vector space instead, i.\,e., \[ \langle gx,y \rangle = \sum_i
      (A^{-1}e'_i)^t X^t G^{-1} Y \star e_i'\] and \[\langle x,g^{-1}y \rangle =
      \sum_i e_i'^t X^t G^{-1} Y A \star e_i'\] for appropriate matrices
    $X,G^{-1},Y,\star$ and $(e_i')_i$ the canonical basis of $K^{d\choose k}.$
    (The matrix $G^{-1}$ will not be invertible in general, even though the
    notation does suggest this.) As $A\star = \star (A^{-1})^t$ by
    \ref{prop:hodge-star-adjoint,prop:compact-group-trivial-lie-twist}, equality
    follows, as the trace is invariant under cyclic permutations.
\end{proof}
\end{theorem}

\section{Tamme's Comparison Results}

We will summarise the results from \cite{MR3352825} which we need as follows:

\begin{theorem}
  \label{prop:tamme-summary}
  \label{prop:analytic-cochains-to-koszul}
  Let $K$ be a complete non-archimedean field of characteristic zero.
  Let $G$ be a Lie group over $K$ and $V$ a barrelled LCVS with an analytic
  action of $G$. Then there is a functorial morphism \[ C^\bullet(G,V) \ra
    \CE^\bullet(\mathfrak g,V)\]  from the analytic cochains
  of $G$ with coefficients in $V$ to the Chevalley-Eilenberg complex
  of the Lie algebra $\mathfrak g$ of $G$ with coefficients in $V$.

  For an open subgroup $U\leq G,$ we denote its Lie algebra by $\mathfrak g(U).$
  Above morphism induces for all $n$ isomorphisms \[\varinjlim_{U\leq_o G}
    H^n(U,V) \cong \varinjlim_U H^n(\mathfrak g(U),V) = H^n(\mathfrak g,V).\]

  The adjoint action of $G$ on $\mathfrak g$ together with the action of $G$ on
  $V$ induce an action of $G$ on the Chevalley-Eilenberg complex and on the Lie
  algebra cohomology groups. If $G$ is compact, then above morphism of complexes
  induces an isomorphism \[ H^n(G,V)\cong H^n(\mathfrak g,V)^G\] for all $n$.
  \begin{proof}
    \cite[sections 3-5]{MR3352825}
  \end{proof}
\end{theorem}

\section{The Duality Theorem}
\begin{lemma}
  \label{prop:finite-grp-inv-eq-coinv}
  Let $G$ be a finite group acting linearly on an $L$-vector space $V$. If the order
  of $G$ is invertible in $L$, the composition of the canonical inclusion and
  projection \[ V^G \ra V \ra V_G\] is an isomorphism.
  \begin{proof}
    By Maschke's theorem, $L[G]$ is a semisimple ring. Therefore there exists an
    $L[G]$-submodule $W$ of $V$ with $V=V^G\oplus W.$ Without loss of generality
    we can assume that $W$ is irreducible. Denote by $I$ the augmentation ideal
    in $L[G]$. Then \[ V_G = V^G/IV^G \oplus W/IW = V^G \oplus W/IW\] and as $W$
    is irreducible, $IW$ is either $0$ or $W$. If $IW=0,$ then $W\subseteq V^G$
    and hence $W=0$ by assumption, so $W/IW=0$ in any case.
  \end{proof}
\end{lemma}

Fix now a complete non-archimedean field $K$ of characteristic zero and a Lie
group $G$ over $K$, which acts equi-analytically on an LCVS $V$.

\begin{lemma}
  \label{prop:finite-generation-of-lie-alg-coh-yields-trivializing-open-subgrp}
  Let $R$ be a $K$-algebra. Assume that $V$ carries the structure of an
  $R$-module, such that the operation of $G$ on $V$ is $R$-linear. If
  $H^i(\mathfrak g,V)$ is finitely generated over $R$, then there is an open
  subgroup of $G$ which acts trivially on $H^i(\mathfrak g,V)$.
  \begin{proof}
    By \ref{prop:tamme-summary}, \[ \varinjlim_{U\leq_o G,\mathrm{res}} H^i(U,V)
      = H^i(\mathfrak g,V),\] which is $R$-linear by our assumptions. Taking
    preimages of the finitely many generators in $H^i(\mathfrak g,V),$ we see
    that there is an open subgroup $U\leq G$ such that $H^i(U,V)\ra
    H^i(\mathfrak g,V)$ is surjective. This $U$ then operates trivially on
    $H^i(\mathfrak g,V)$.
  \end{proof}
\end{lemma}

\begin{theorem}
  \label{prop:duality-morphism-analytic-cohomology}
  If $G$ is compact and $V,V'_b$ barrelled, we get a functorial (in $V$)
  morphism of complexes \[ C^\bullet(G,V'_b) \ra \Hom_K(C^\bullet(G,V),
    K)[-d].\]

  If one of the following two
  conditions is satisfied:
  \begin{itemize}
  \item An open subgroup of $G$ operates trivially on the Lie algebra
    cohomology, the differentials in the Chevalley-Eilenberg complex are strict
    and Hahn-Banach holds for $V$, or
  \item $V$ is finite-dimensional,
  \end{itemize}
  then this morphism induces
  isomorphisms \[ H^i(G,V'_b)\cong H^{d-i}(G,V)'.\]
  \begin{proof}
    Note first that by
    \ref{prop:finite-generation-of-lie-alg-coh-yields-trivializing-open-subgrp},
    an open subgroup of $G$ operates trivially on the Lie algebra cohomology, no
    matter the case.
    
    By
    \ref{prop:analytic-cochains-to-koszul}
    we have morphisms \[ C^\bullet(G,V'_b)\ra \CE^\bullet(\mathfrak g,V'_b)\]
    and \[ \Hom_K(\CE^\bullet(\mathfrak g,V), K) \ra \Hom_K(C^\bullet(G,V),K).\]

    As $G$ is compact, we can employ \ref{prop:chevalley-eilenberg-duality} without
    having to twist the Lie algebra action (cf.~\ref{prop:compact-group-trivial-lie-twist}).
    We therefore get a $G$-equivariant isomorphism \[ \CE^\bullet(\mathfrak g, V'_b) \cong
      \Hom_{K,\mathrm{cts}}(\CE^{\bullet}(\mathfrak g, V), K)[-d].\] Composition
    with the inclusion
    \[ \Hom_{K,\mathrm{cts}}(\CE^{\bullet}(\mathfrak g, V), K)\subseteq
      \Hom_{K}(\CE^{\bullet}(\mathfrak g, V), K)\] then yields the comparison
    morphism, which is clearly functorial in $V$. 
    If the differentials in the complex $\CE(\mathfrak g, V)$ are strict and
    Hahn-Banach holds for $V$, we get a $G$-equivariant isomorphism on the level
    of cohomology: \[ H^i(\mathfrak g,V'_b) \ra
      \Hom_{K,\mathrm{cts}}(H^{d-i}(\mathfrak g, V), K).\] Especially we get the
    following commutative diagram:
    \[
      \begin{tikzcd}
                                     &  &
                                     &  & {(H^{d-i}(\mathfrak g, V)^G)'}
                                     \arrow[rrd, "\cong", dashed]      &  &                 \\
{H^i(G,V'_b)} \arrow[rr] \arrow[rrd, "\cong", dashed] & & {H^i(\mathfrak
  g,V'_b)} \arrow[rr, "\cong"]             &  & {H^{d-i}(\mathfrak g, V)'} \arrow[u, two heads] \arrow[rr] &  & {H^{d-i}(G,V)'} \\
                                     &  & {H^i(\mathfrak g,V'_b)^G} \arrow[rr,"\cong"]
                                     \arrow[u, hook] &  & {(H^{d-i}(\mathfrak g, V)_G)'} \arrow[u]        &  &                
\end{tikzcd}
    \]
    The dashed isomorphisms are again instances of
    \ref{prop:analytic-cochains-to-koszul}.
    If an open subgroup of $G$ acts trivially on the Lie algebra cohomology,
    then the composition \[ (H^{d-i}(\mathfrak g, V)_G)'\ra H^{d-i}(\mathfrak g,
      V)\ra (H^{d-i}(\mathfrak g, V)^G)' \] is an isomorphism by
    \ref{prop:finite-grp-inv-eq-coinv} and the claim follows.
  \end{proof}
\end{theorem}

\begin{corollary}
  \label{prop:analytic-duality-fin-dim}
  Let $G$ be a compact Lie group of dimension $d$ acting analytically on a
  finite dimensional $K$-vector space $V$. Then we have a functorial
  quasi-isomorphism
  \[ C^\bullet(G,V^*) \oset{\cong}{\ra} C^\bullet(G,V)^*[-d].\]
  \begin{proof}
    By \ref{prop:fin-dim-is-equianalytic}, we are in the setting of
    \ref{prop:duality-morphism-analytic-cohomology}. If $V$ is
    finite-dimensional, we see that $\CE(\mathfrak g, V')$ is a complex of
    finite dimensional vector spaces and the analytic cohomology groups are
    therefore finite-dimensional as well by \ref{prop:tamme-summary}. For all
    cohomology groups involved, their abstract duals hence coincide with their
    continuous duals and the result follows.
  \end{proof}
\end{corollary}

\begin{remark}
  \label{rmk:functoriality-of-duality-morphism}
  Functoriality in \ref{prop:duality-morphism-analytic-cohomology} means the
  following: Let $V,W$ be LCVS with equi-analytic actions of $G$ on
  them. Assume that $V,W,V'_b,W'_b$ are all barrelled. Given a $G$-equivariant
  continuous linear map $\varphi\colon V\ra W,$ we get a commutative diagram:
  \[
    \begin{tikzcd}[column sep=large]
      C^\bullet(G,W'_b) \arrow{r}\arrow{d}{C(G,\varphi')} & \Hom_K(C^\bullet(G,W),K)[-d]\arrow{d}{\Hom_K(C(G,\varphi),K)}\\
      C^\bullet(G,V'_b) \arrow{r} & \Hom_K(C^\bullet(G,V),K)[-d]
    \end{tikzcd}
  \]
  That the maps involved are well-defined follows from
  \ref{prop:dual-is-functor,prop:cts-after-analytic-is-analytic}.
\end{remark}

\begin{remark}
  Of course a quasi-isomorphism would be a nicer result in the setting of
  \ref{prop:duality-morphism-analytic-cohomology}.
  The obvious strategy would be topologise $C^\bullet(G,V)$ in such a manner
  that the differentials are strict and that the cohomology groups are
  topologically identical to the topology on the Lie algebra cohomology. The
  same argument as above would then (under the additional hypotheses on the
  Chevalley-Eilenberg complex and the Lie algebra cohomology) yield a
  quasi-isomorphism
  \[ C^\bullet(G,V')\ra\Hom_{K,\mathrm{cts}}(C^\bullet(G,V),K).\] We consider
  this endeavour to be rather futile, which is the reason we axiomatised and
  capsuled topological considerations in
  \cite{thomas:cohomology-of-topologised-monoids} in the first place.
\end{remark}
\begin{remark}
  \label{rmk:duality-is-poincare-duality}
  If $K=\Q_p$ and $V$ is a finite-dimensional $\Q_p$-vector space, then by
  \cite[V.(2.3.10)]{MR0209286} analytic cohomology is just continuous
  cohomology. \Cref{prop:duality-morphism-analytic-cohomology} is then a
  possible way to phrase Poincaré duality, which however does \emph{not}
  coincide with Poincaré duality due to Lazard
  (cf.~\cite[V.(2.5.8)]{MR0209286}). Poincaré duality there is an
  \emph{integral} phenomenon and the dual is given by
  $\Hom_{\Z_p,\mathrm{cts}}(V,\Q_p/\Z_p)$.
\end{remark}

\begin{example}
  Let $V$ be any barrelled LCVS and $G$ a compact abelian Lie group over $K$ of
  dimension $d$. The trivial action of $G$ on $V$ is of course
  equi-analytic. The Lie algebra $\mathfrak g$ of $G$ then operates by zero on
  $V.$ The differentials in the Chevalley-Eilenberg complex
  $\CE^\bullet(\mathfrak g,V)$ are all zero. \Cref{prop:tamme-summary} then
  yields that $H^i(G,V)\cong H^i(\mathfrak g,V) = \Hom_K(\bigwedge^i\mathfrak
  g,V)$ and the isomorphism $H^i(G,V'_b)\cong H^{d-i}(G,V)'$ of
  \ref{prop:duality-morphism-analytic-cohomology} stems from the pairing \[
    \wedge\colon\textstyle \bigwedge^i\mathfrak g \times \bigwedge^{d-i}\mathfrak g \ra
    K. \]
\end{example}

\section{Two Applications to $(\varphi,\Gamma)$-modules}

Analytic cohomology as it appears in the theory of $(\varphi,\Gamma)$-modules
has mostly had a very strong \emph{ad hoc} flavour. Arguments often used
crucially that $\Gamma$ is a one-dimensional Lie group and $\varphi$ a single
operator. However, our framework of
\cite{thomas:cohomology-of-topologised-monoids} and our results of the previous
section are much more flexible and easily apply themselves to higher-dimensional
$\Gamma$ and multiple operators $\varphi_1,\dots,\varphi_d.$ The natural objects
to look at are thus \emph{multivariable} $(\varphi,\Gamma)$-modules. There is,
however, no unified notion of multivariable $(\varphi,\Gamma)$-modules and to
our knowledge, no definition of multivariable Lubin-Tate
$(\varphi,\Gamma)$-modules has been published. Consequently many results which
are well known in the univariate case are unknown to hold in the multivariable
case. Our arguments, not relying on the ad hoc constructions of analytic
cohomology, should easily carry over to the multivariable case as soon as the
necessary category equivalences are shown. An important step towards this has
recently been done by \cite{pal-zabradi:cohomology-and-overconvergence}.

\subsection{An Exact Sequence of Berger and Fourquaux}

We start by improving a result of Berger and Fourquaux, cf.~\cite[theorem
2.2.4]{MR3690276}, which can be stated without precisely defining
$(\varphi,\Gamma)$-modules.

Let $F|\Q_p$ be a local field and consider the category $\cat{M}$ of analytic
$F$-manifolds.

We fix an LF-space $ A \cong \varinjlim_r\varprojlim_s A^{[r,s]}$ with Banach
spaces $A^{[r,s]}$ for the remainder of this subsection (cf.~\ref{def:lf-space}).
The notation $A^{[r,s]}$ will become apparent in the next subsection.

\begin{definition}
  \label{prop:pro-analytic-direct-limit}
  For $X\in\cat{M}$ let $f\colon X\ra A$ be a continuous map. We call $f$
  pro-$F$-analytic, if there exists an $r$ and a factorisation \[f\colon X\ra
    \varprojlim_s A^{[r,s]}\ra A,\] such that all induced maps \[\begin{tikzcd}
      X\arrow{r} & \varprojlim A^{[r,s]}\arrow{d}\\ & A^{[r,s']} \end{tikzcd}\] are
  locally $F$-analytic. We denote the set of pro-$F$-analytic maps from $X$ to
  $A$ by $h(X,A)$, i.\,e., \[ h(X,A) = \varinjlim_r\varprojlim_s
    h_{\mathrm{an}}(X,A^{[r,s]}),\] where $h_{\mathrm{an}}$ denote the analytic
  maps in the sense of \ref{def:analytic-map}.
\end{definition}
\begin{proposition}
  \label{prop:analytic-is-pro-analytic}
  An $F$-analytic map into a Fréchet space in the
  sense of \ref{def:analytic-map} is pro-$F$-analytic. 
  \begin{proof}
    Let $B=\varprojlim_n B_n$ be a Fréchet space with all $B_n$ Banach.
    We need to show that if $f\colon M\ra B$ is analytic, then so is $f\colon
    M\ra B\ra B_n$ for each $n$. But this is precisely the content of \ref{prop:cts-after-analytic-is-analytic}.
  \end{proof}
\end{proposition}
\begin{remark}
  The argument in \ref{prop:cts-after-analytic-is-analytic} shows why not every
  pro-analytic map needs to be analytic: For a pro-analytic map
  $M\ra\varprojlim_n B_n$ (and around a fixed point in $M$), we have a positive
  radius of convergence $R_n$ of the power series development for every $B_n$.
  But there is no need for $\inf_n R_n$ to be positive, which is the natural
  estimate for the radius of convergence for $B.$
\end{remark}

Let $\Gamma$ be an analytic group over $F$ and $A$ an LF-space an action from
$\Gamma$ by pro-$F$-analytic maps and a continuous $L$-linear endomorphism
(which we will call $\psi$) of $A$. $A$ is then a $G=\psi^{\N_0}\times\Gamma$-module.
By setting $C^n(G,A)=h(G^n,A)$ and using the usual inhomogeneous cochain
differential, we can define the pro-$F$-analytic cohomology of $G$ with
coefficients in $A$ as the cohomology of this complex. This yields a
well-defined cohomology theory that exhibits many of the standard features,
cf.~\cite{thomas:cohomology-of-topologised-monoids}.

We now prove the following stronger version of \cite[theorem 2.2.4]{MR3690276},
where only the case of one-dimensional $\Gamma$ is considered. For the
one-dimensional case, we can also show that the last map appearing in the exact
sequence of Berger and Fourquaux is surjective.

\begin{theorem}
  \label{prop:ex-seq-bf}
  There is an exact sequence
  of pro-$F$-analytic cohomology groups:
  \[
    \begin{tikzcd}
      0 \arrow{r} & H^1(\Gamma, A^{\psi=1}) \arrow{r} & H^1(\psi^{\N_0}\times\Gamma, A)\arrow{r} &
      (A/(\psi-1)A)^\Gamma  \arrow[looseness=1.2, overlay, out=355, in=175]{dll}\\
      & H^2(\Gamma, A^{\psi=1}) \arrow{r}& H^2(\psi^{\N_0}\times\Gamma, A)   &
    \end{tikzcd}
  \] The second to last
  group can be replaced by
  zero if $\Gamma$ is compact, has dimension one over $L$, and also operates
  analytically on each $\varprojlim_sA^{[r,s]}.$
  \begin{proof}
    In lieu of \cite[theorem 10.26]{thomas:cohomology-of-topologised-monoids} it
    only remains to show the surjectivity onto $(A/(\psi-1)A)^\Gamma$ for
    compact one-dimensional $\Gamma$. Note that we can't directly use
    \ref{prop:tamme-summary} to compare this cohomology group to Lie algebra
    cohomology, as pro-analytic maps and analytic maps don't necessarily
    coincide.

    By \ref{prop:pro-analytic-direct-limit} and the exactness of direct limits
    we can assume that $A$ is Fréchet. As every analytic map is also
    pro-analytic, the same proof as for \cite[proposition
    7.1]{thomas:cohomology-of-topologised-monoids} together with \cite[theorem
    10.26]{thomas:cohomology-of-topologised-monoids} shows that we have the
    following commutative diagram with exact rows:
    \[
\begin{tikzcd}[column sep=5ex]
  0 \arrow[r] & {H^1(\Gamma, A^{\psi=1})} \arrow[r]                & {H^1(\psi^{\N_0}\times\Gamma, A)} \arrow[r]                &  (A/(\psi-1)A)^\Gamma \arrow[r]                 & {H^2(\Gamma, A^{\psi=1})}          \\
  0 \arrow[r] & {H^1_{\mathrm{an}}(\Gamma, A^{\psi=1})} \arrow[r] \arrow[u,
  hook] & {H^1_{\mathrm{an}}(\psi^{\N_0}\times\Gamma, A)} \arrow[r] \arrow[u,
  hook] & (A/(\psi-1)A)^\Gamma \arrow[r] \arrow[u, equal] &
  {H^2_{\mathrm{an}}(\Gamma, A^{\psi=1})} \arrow[u]
\end{tikzcd}
    \] Here, $H^\bullet_{\mathrm{an}}$ denotes cohomology with respect to
    analytic maps in the sense of \ref{def:analytic-map}, for which the results
    of \cite{thomas:cohomology-of-topologised-monoids} also hold.
    For these cohomology groups, we can apply \ref{prop:tamme-summary} to see
    that \[{H^2_{\mathrm{an}}(\Gamma, A^{\psi=1})}=0. \] The exact sequence
    \[0 \ra H^1(\Gamma, A^{\psi=1}) \ra H^1(\psi^{\N_0}\times\Gamma, A)\ra
      (A/(\psi-1)A)^\Gamma \ra 0\] follows at once. 
  \end{proof}
\end{theorem}

\begin{remark}
  \label{rmk:we-hate-lim1}
  Assume that $\Gamma$ is compact and of dimension one and that $A=\varprojlim_n
  A_n$ is Fréchet. Then even if the operation on $A$ is only pro-analytic, the
  same argument as in \ref{prop:ex-seq-bf} yields exact
  sequences \[\begin{tikzcd} 0\arrow{r} & H^1_{\mathrm{an}}(\Gamma,
      A_n^{\psi=1}) \arrow{r} & H^1_{\mathrm{an}}(\psi^{\N_0}\times\Gamma,
      A_n)\arrow{r} & (A_n/(\psi-1) A_n) \arrow{r} & 0\end{tikzcd}\] for every
  $n$. Assume furthermore that the image of $A_m$ in $A_n$ is dense for every
  $m\geq n$. We then have isomorphisms \[ H^1(\psi^{\N_0}\times\Gamma, A) \cong
    \varprojlim_n H^1_{\mathrm{an}}(\psi^{\N_0}\times\Gamma, A_n) \] by
  \cite[proposition 2.1.1]{MR3690276}. As taking invariants commutes with
  projective limits, we therefore get the following exact
  sequence: \[\begin{tikzcd}[column sep=4ex] 0 \arrow{r} & H^1(\Gamma,
      A^{\psi=1}) \arrow{r} & H^1(\psi^{\N_0}\times\Gamma,A) \arrow{r} &
      (A/(\psi-1)A)^\Gamma \arrow{r} & \varprojlim_n\!\!\!\!{}^1\,
      H^1_{\mathrm{an}}(\Gamma, A_n^{\psi=1})\end{tikzcd}\] Using again
  \ref{prop:tamme-summary}, we see that $H^1_{\mathrm{an}}(\Gamma,
  A_n^{\psi=1})\cong H^1(\mathfrak g, A_n^{\psi=1})^\Gamma,$ where $\mathfrak g$
  is the Lie algebra of $\Gamma.$ The action of $\Gamma$ on $\mathfrak g$ is
  trivial and $\mathfrak g\cong L,$ so $ H^1(\mathfrak g, A_n^{\psi=1})^\Gamma$
  has a comparatively simple description as the $\Gamma$-invariants of certain
  quotients of $A_n^{\psi=1}$, which depend on the precise group action,
  cf.~e.\,g.~\cite[theorem 7.4.7]{MR1269324}. For these it might in certain
  examples be possible to show the (topological) Mittag-Leffler condition,
  cf.~e.\,g.~\cite[(0.13.2.4)]{MR217085}, and hence show that \[0 \ra
    H^1(\Gamma, A^{\psi=1}) \ra H^1(\psi^{\N_0}\times\Gamma, A)\ra
    (A/(\psi-1)A)^\Gamma \ra 0\] is exact.
\end{remark}

\begin{remark}
  \label{rmk:proanalytic-vs-analytic-cohomology}
  Considering \ref{rmk:we-hate-lim1}, it is natural to ask for the relationship
  between $ H^k(\Gamma, A)$ and $\varprojlim_n H^k(\Gamma, A_n)$, where
  $A=\varprojlim_n A_n$ is again assumed to be Fréchet.

  Consider the exact sequence
    \[\begin{tikzcd}[column sep=4ex]
        0\arrow{r} & C^\bullet(\Gamma, A) \arrow{r} & \prod_n
        C^\bullet_{\mathrm{an}}(\Gamma, A_n) \arrow{r}{1-u} & \prod_n
        C^\bullet_{\mathrm{an}}(\Gamma,A_n) \arrow{r} &
        \varprojlim_n\!\!\!\!{}^1\, C^\bullet_{\mathrm{an}}(\Gamma,A_n)
        \arrow{r} & 0,\end{tikzcd}\] whose middle map is given by
    \[
      1-u\colon (f_n)_n\mapsto (f_n-(A_{n+1}\shortrightarrow A_n)\circ f_{n+1})_n.
    \]
    Its existence follows from very general arguments,
    cf.~\cite[(2.4.7)]{MR2392026} for a correct statement and proof.
  
    If one could show that $1-u$ is indeed surjective, then the long exact
    sequence of cohomology would yield the following short exact sequences for
    every $k$:
    \[
      \begin{tikzcd}
        0 \arrow{r} &\varprojlim_n\!\!\!\!{}^1\, H^{k-1}_{\mathrm{an}}(\Gamma,
        A_n)\arrow{r} & H^k(\Gamma, A) \arrow{r} & \varprojlim_n
        H^k_{\mathrm{an}} (\Gamma, A_n) \arrow{r} & 0
      \end{tikzcd}
    \]
    Write $d=\dim_F\Gamma.$ Then these short exact sequences would show that
    $H^k(\Gamma,A)=0$ for every $k>d+1$ and that $H^{d+1}(\Gamma,A)\cong
    \varprojlim_n\!\!\!\!{}^1\,H^d_{\mathrm{an}}(\Gamma,A_n).$ Especially, the
    considerations in \ref{prop:ex-seq-bf,rmk:we-hate-lim1} would coincide for
    $d=1.$
\end{remark}


\subsection{So Many Rings}
\label{sec:so-many-rings}

Fix a complete field $L\subseteq \C_p.$ We mostly follow the notation of \cite{MR3690276}. 
\begin{definition}
  Consider the abelian group $\Map(\Z, L) = L[[X, X^{-1}]]$. In this set, we can
  define the following rings. Let $I\subset [0,1]$ be an interval. Set
  \begin{gather*}
    \Brs{L}{I} = \left\{ \sum_{i\in\Z} a_{i}
      X^{i} \in L[[X,X^{-1}]]\mid \text{convergent for } z\in\C_p, \abs{z}\in I
    \right\}.
  \end{gather*}
  For $I=[r,s],$ this is a Banach space over $L$ with norm
  $\norm{-}_{\Brs{L}{[r,s]}}=\max \{\norm{-}_r,\norm{-}_s\},$ where $\norm{-}_r$
  and $\norm{-}_s$ are defined via \[ \norm{\sum_{i\in\Z}a_i X^i}_t =
    \sup_{i\in\Z} \abs{a_i}t^i.\] We also define
  \begin{gather*}
    \Brigdeath{L}{r} = \varprojlim_{r<s<1,s\rightarrow 1}\Brs{L}{[r,s]} \text{ for } 0<r<1,
    \text{ and}\\
    \Brigdeath{L} = \varinjlim_{0<r<1, r\rightarrow 1}\Brigdeath{L}{r}.
  \end{gather*}
  For $L$ finite over $\Q_p$, we can also define the complete discrete valuation ring
  \begin{gather*}
    \Ars{L} =\mathcal O_L
    [[X]][X^{-1}]^{\wedge}\end{gather*} with quotient field \begin{gather*}
    \Brs{L} = \Ars{L}[p^{-1}],
  \end{gather*}
  where $-^\wedge$ denotes $p$-adic completion. We also define \[\Bdeath{L} =
    \left\{ f\in\Brs{L} \mid f \text{ has a non-empty domain of convergence}
    \right\}.\]
  If $M|L$ is a finite extension, the theory of the field of norms provides a
  certain ring extension $\Ars{M|L}$ over $\Ars{L}.$ Its quotient field
  $\Brs{M|L}$ is an unramified extension of $\Brs{L}$. We define the following complete
  discrete valuation ring and its quotient field:
  \begin{gather*}
    \Ars = (\bigcup_{M|L\text{ finite}}\Ars{M|L})^\wedge,\\
    \Brs = (\bigcup_{M|L\text{ finite}}\Brs{M|L})^\wedge.
  \end{gather*}
\end{definition}

\begin{remark}
  Some authors denote the rings $\Ars{M|L}$ and $\Brs{M|L}$ simply by $\Ars{M}$
  and $\Brs{M}$, obfuscating the fact that these rings are \emph{relative}
  notions: In our notation, we always have $\Ars{M|M}=\Ars{M}$ but generally
  $\Ars{M|L}\neq \Ars{M}.$ We consider this abuse of notation in the literature
  truly abusive.
\end{remark}

\subsection{Lubin-Tate $(\varphi,\Gamma)$-modules}
\label{sec:lubin-tate-phi-gamma-modules}
Fix a finite Galois extension $L|\Q_p$ for the remainder of this chapter. We
denote the ring of integers of $L$ by $\mathcal O_L$ and its residue field of
cardinality $q$ by
$\kappa_L$. We also fix a uniformiser $\pi$ of $\mathcal O_L$.

\subsubsection{The Lubin-Tate Case}
We assume familiarity with the theory of formal multiplication in local fields, 
cf.~e.\,g.~\cite[section 3]{MR0220701}.

Denote by $LT$ the Lubin-Tate formal $\mathcal O_L$-module attached to $\pi$,
i.\,e., as a set $LT$ is the maximal ideal of the integral closure of $\mathcal
O_L$ in an algebraic closure of $L$ and the addition is defined via the unique
formal group law corresponding to the endomorphism \[ [\pi](T)=T^q+\pi T.\]
Lubin-Tate theory then yields commuting power series $[a](T)\in\mathcal O[[T]]$
for all $a\in\mathcal O_L$, which give rise to a $\mathcal O_L$-module structure
on $LT.$ It also yields a homomorphism \[ \chi_{LT}\colon G_L\ra \mathcal
  O_L^\times,\] which induces an isomorphism \[\chi_{LT}\colon
  \Gamma_L=G(L_\infty|L)\oset{\cong}{\ra} \mathcal O_L^\times,\] where
$L_\infty$ is the extension of $L$ generated by all $\pi^\infty$-torsion points
of $LT$.

For $f(T)$ in any of the rings $\Brigdeath{L},\Bdeath{L},\Ars{L},\Brs{L}$ we
have well-defined elements
\begin{gather*} \varphi(f)(T) = f([\pi](T)),\\
  (gf)(T) = f([\chi_{LT}(g)](T)),\, g\in\Gamma_L.
\end{gather*}

Denote the monoid $\varphi^{\N^0}$ by $\Phi.$ Then by construction above formula
induce a continuous action of $\Phi\times\Gamma_L$ by ring homomorphisms
on each of the above rings with their respective topologies, which for
$\Brigdeath{L}$ is even pro-$L$-analytic
cf.~e.\,g.~\cite[theorem 8.1]{MR3577371}.

\begin{definition}
  Let $R$ be either of $\Brigdeath{L},\Bdeath{L},\Ars{L},\Brs{L}$. A
  $(\varphi,\Gamma_L)$-module $M$ over $R$ is a free $R$-module of finite rank
  with a semi-linear continuous action of $\Phi\times\Gamma_L.$ It is called étale
  if $\varphi(M)$ generates $M$.
\end{definition}

\subsubsection{Relation to Iwasawa Cohomology}

Recall the following result due to Kisin and Ren.

\begin{theorem}[{\cite[theorem 1.6]{MR2565906}}]
  The functor \[ V\mapsto \D_{\mathcal O}(V) = (\Ars\otimes_{\mathcal O_L}
    V)^{G(\overline{\Q}_p|L_\infty)}\] establishes an equivalence between the
  categories of $\mathcal O_L$-linear representations of $G_L$ and étale
  $(\varphi,\Gamma_L)$-modules over $\Ars{L}.$
\end{theorem}

For any étale $(\varphi,\Gamma_L)$-module $D$ over $\Ars{L}$ there is an $\mathcal
O_L$-linear endomorphism \[ \psi\colon D\ra D\] satisfying \[ \psi\circ\varphi =
  \frac{q}{\pi}\id_D,\] cf.~e.\,g.\,\cite[416]{MR3629658}.

There is the following relationship between $(\varphi,\Gamma_L)$-modules and
Iwasawa cohomology. While $(\varphi,\Gamma_L)$-modules have plentiful
applications, this is our main reason for studying them.

\begin{theorem}[{\cite[theorem 5.13]{MR3629658}}]
  \[ H^1_{\mathrm{Iw}}(L_\infty, V(\chi_{LT}^{-1}\chi_{cyc})) = \D_{\mathcal
      O}(V)^{\psi=1},\]
  where $\chi_{cyc}$ denotes the cyclotomic character.
\end{theorem}

\begin{corollary}
  If $L\neq\Q_p$, we have $D^{\psi=1,\Gamma=1}=0$ for any
  étale $(\varphi,\Gamma_L)$-module $D$ over $\Ars{L}$. 
  \begin{proof}
    Together with these two aforementioned results, this follows immediately
    from
    \cite[theorem 8.2]{thomas-venjakob:on-spectral-sequences-for-iwasawa-adjoints-a-la-jannsen-for-families},
    as elements in $D^{\psi=1,\Gamma=1}$ are torsion over the Iwasawa algebra. 
  \end{proof}
\end{corollary}

\subsubsection{Overconvergence}

For the remainder of this chapter, we also fix a finite extension $F|L.$

\begin{definition}
  For an $L$-linear representation $V$ of $G_F$ set \[
    \mathrm{D}(V) = \left( \Bdeath\otimes_L V
    \right)^{G(\overline{\Q}_p|FL_\infty)}.\] $\Gamma = G(FL_\infty|F)$ is an
  open subgroup of $\Gamma_L$ and by the Lubin-Tate character hence isomorphic to
  an open subgroup of $\mathcal O_L^\times.$ $\Drs(V)$ is an étale
  $(\varphi,\Gamma)$-module over $\Brs{L}$.
\end{definition}
\begin{definition}
  Let $D$ be a $(\varphi,\Gamma)$-module over $\Brs{L}.$ If there is a basis of
  $D$ such that all endomorphisms in $\Phi\times\Gamma$ have representation
  matrices in $\Bdeath{L},$ we call $D$ overconvergent. This basis generates a
  $(\varphi,\Gamma)$-module over $\Bdeath{L}$, which we will call $\Ddeath$. A
  Galois representation $V$ is called overconvergent if $\D(V)$ is. We
  will then write $\Ddeath(V)$ instead of $\D(V)^\dagger.$
\end{definition}

\begin{definition}
  Let $V$ be an overconvergent Galois representation. Set \[ \Drigdeath(V) =
    \Brigdeath{F}\otimes_{\Bdeath{F}}\Ddeath(V).\]
\end{definition}
\begin{definition}
  A finite dimensional $L$-linear representation $V$ of $G_F$ is called
  $L$-analytic, if for all embeddings $\tau\colon L\ra\bar{\Q}_p$ different from
  the fixed one, \[\C_p\otimes_L^\tau V\] is a trivial semilinear
  $\C_p$-representation, i.\,e., as a Galois module it is isomorphic to
  $\C_p\otimes_L\widetilde{V}$ for an $L$-vector space $\widetilde{V}$ with
  trivial Galois operation.
\end{definition}

\begin{lemma}
  \label{prop:dual-is-analytic}
  A finite dimensional $L$-linear representation $V$ is $L$-analytic if and only
  if $V^*=\Hom_L(V,L)$ is.
  \begin{proof}
    Note that a representation is trivial if and only if its dual is. The
    statement then follows from the isomorphisms \[ (\C_p\otimes^\tau_F V)^\vee
      \cong \C_p^\vee \otimes^\tau_F V^* \cong \C_p\otimes^\tau_F V^*,\] where
    $-^\vee=\Hom_{\C_p}(-,\C_p).$
  \end{proof}
\end{lemma}

\begin{proposition}
  The action of $\Gamma$ on $\Drigdeath(V)$ is pro-$L$ analytic.
  \begin{proof}
    This follows for example from \cite[theorem 8.1]{MR3577371}.
  \end{proof}
\end{proposition}

\begin{proposition}[{\cite[2554]{MR3194651}}]
  Let $D$ be an étale $(\varphi,\Gamma)$-module over $\Brigdeath{F},$ which be
  finite generation can be written as $D=\Brigdeath{F}\otimes_{\Brigdeath{F}{r}}
  D^r$ for some $r$ and some $(\varphi,\Gamma)$-module $D^r$ over
  $\Brigdeath{F}{r}.$ Then the series \[ \log g = \sum_{i=1}^\infty
    \frac{(-1)^{i-1}}{i}(g-1)\] converges for $g$ small enough to an operator on
  $D^r.$ By $\Z_p$-linear extension, this gives rise to a well-defined action of
  the Lie algebra $\mathfrak g\cong L$ via \[(x,d)\mapsto (\log(\exp x))(d).\]
\end{proposition}

\begin{definition}
  A $(\varphi,\Gamma)$-module $D$ over $\Brigdeath{F}$ is called $L$-analytic,
  if the action of the Lie algebra of $\Gamma$ on $D$ is $L$-linear.
\end{definition}

Berger then shows the following refinement of the category equivalence.

\begin{theorem}[{\cite[theorem D]{MR3577371}}]
  \label{prop:equivalence-analytic-phi-gamma-modules}
  $V\mapsto \Drigdeath(V)$ is an equivalence of categories between $L$-analytic
  representations of $G_F$ and étale $L$-analytic $(\varphi,\Gamma)$-modules
  over $\Brigdeath{F}.$
\end{theorem}

\subsection{Towards Duality in the Herr Complex}
\label{sec:herr-duality}
We continue to use the notation from \ref{sec:lubin-tate-phi-gamma-modules}.

Let $V$ be an $L$-analytic representation of $G_F$. Then the Herr complex is
given by the double complex
\[ C^\bullet(\Gamma,\Drigdeath(V)) \oset{\varphi -1}{\ra}
  C^\bullet(\Gamma,\Drigdeath(V)),\] whose attached double complex is
quasi-isomorphic to \[
  C^\bullet(\Phi\times\Gamma,\Drigdeath(V))\] by
\cite[theorem 11.6]{thomas:cohomology-of-topologised-monoids}. Here $C^\bullet$ denotes the
pro-$F$-analytic cochains.

It seems natural to try to apply our duality result
\ref{prop:duality-morphism-analytic-cohomology} to this setting, however, this
is not immediately possible.

Starting with $\Drigdeath(V)$, there are (at least) three natural ways to
dualise this object: We can consider the $(\varphi,\Gamma)$-module attached to
the dual representation, the module theoretic dual over the Robba ring
$\Brigdeath{L}$, or the topological dual $\Drigdeath(V)'_b.$ We first need to
investigate how they relate to one another.

\begin{lemma}
  \label{prop:identify-duals-over-robba}
  $\Drigdeath(V^*) = \Hom_{\Brigdeath{L}}(\Drigdeath(V),\Brigdeath{L}).$
  \begin{proof}
    It is well known that the category equivalence over $\Brs{L}$ is a functor
    of closed monoidal categories and hence respects taking duals. As the duals
    of analytic representations are again analytic by
    \ref{prop:dual-is-analytic}, the finer category equivalence
    \ref{prop:equivalence-analytic-phi-gamma-modules} also has to respect duals.
  \end{proof}
\end{lemma}

Let $\Omega\in\C_p$ be the period of $LT,$ cf.~\cite[§1.1.3]{MR3522263}. If
$L\neq\Q_p,$ it is transcendent over $\Q_p.$ Let $K$ be the complete subfield of
$\C_p$ generated by $L$ and $\Omega.$
For a $(\varphi,\Gamma)$-module $D$ over $\Brigdeath{L}$ write $D_K$ for the
respective $(\varphi,\Gamma)$-module over $\Brigdeath{K}$ after extension of
scalars.

Serre duality implies the following result:

\begin{proposition}
  \label{prop:robba-dual-is-top-dual}
 $\Hom_{\Brigdeath{K}}(\Drigdeath(V)_K,\Brigdeath{K})(\chi_{LT})=(\Drigdeath(V)_K)'_b,$
 where the dual is taken over $K$.
 \begin{proof}
   \cite[lemma 2.37]{schneider-venjakob:regulator-maps}
 \end{proof}
\end{proposition}

Note that as the extension $K|\Q_p$ is infinite, we cannot assume that $K$ is
spherically complete. However, we at least have the
following.

\begin{proposition}
  \label{prop:phi-gamma-modules-are-of-countable-type}
  Étale $(\varphi,\Gamma)$-modules over $\Brigdeath{K}$ are of countable type.
  \begin{proof}
    It suffices to show that $\Brigdeath{K}$ is of countable type.
    The completions of $\Brigdeath{K}$ at the various continuous seminorms are
    exactly the rings $\Brigdeath{K}{[r,s]}.$ We will show that the countable
    set $\widetilde{L}(\Omega)[X,X^{-1}]$ is dense in every
    $\Brigdeath{K}{[r,s]}$, where $\widetilde{L}$ is a number field which is
    dense in $L$.

    Let $f=\sum_n a_n X^n\in\Brigdeath{K}{[r,s]}$ and $\varepsilon>0.$
    Convergence of $f$ on the closed annulus of inner radius $r$ and outer
    radius $s$ implies \[ \sup_{n<k} \abs{a_n} r^n \ra 0 \quad (k\to
      -\infty)\] and \[ \sup_{n>k} \abs{a_n} s^n \ra 0 \quad (k\to
      \infty).\] We can therefore choose $k$ with
    $\norm{\sum_{n<-k} a_n X^n}_r, \norm{\sum_{n>k} a_n X^n}_s<\varepsilon$. As
    $\widetilde{L}(\Omega)$ is dense in $K$, we can also choose
    $\alpha_{-k},\dots,\alpha_k\in\widetilde{L}(\Omega)$ such that \[ \max_i
      \abs{a_i-\alpha_i}< \varepsilon r^k.\] It follows that \[
      \norm{f-\sum_{i=-k}^k \alpha_i X^i}_{\Brigdeath{K}{[r,s]}}<\varepsilon.\]
  \end{proof}
\end{proposition}
\begin{proposition}
  \label{prop:partial-herr-duality}
  There are natural morphisms of complexes \[ C^\bullet(\Psi\times\Gamma,
    \Drigdeath(V^*)_K(\chi_{LT})) \ra
    \Hom_K(C^\bullet(\Phi\times\Gamma,\Drigdeath(V)_K),K)[-2]\] and \[
    C^\bullet(\Phi\times\Gamma, \Drigdeath(V^*)_K(\chi_{LT})) \ra
    \Hom_K(C^\bullet(\Psi\times\Gamma,\Drigdeath(V)_K),K)[-2]\] stemming from a
  comparison of Lie algebra cohomology.
  \begin{proof}
    $\Drigdeath(V^*)_K(\chi_{LT})$ is an étale $(\varphi,\Gamma)$-module over
    $\Brigdeath{K}$ and has the structure of an LF-space over $K$, so
    \[ \Drigdeath(V^*)_K (\chi_{LT})= \varinjlim_r\varprojlim_s D^{*,[r,s]},\] where
    $D^{*,[r,s]}$ are Banach spaces over $K$.

    We see that \[C^\bullet(\Gamma,
      \Drigdeath(V^*)_K(\chi_{LT}))=\varinjlim_r\varprojlim_s
      C^\bullet_{\mathrm{an}}(\Gamma,D^{*,[r,s]}).\] By
    \ref{prop:tamme-summary} we get a morphism \[ C^\bullet(\Gamma,
      \Drigdeath(V^*)_K(\chi_{LT}))\ra \varinjlim_r\varprojlim_s \CE^\bullet(\mathfrak g,
      D^{*,[r,s]}).\] Now
    \[\varinjlim_r\varprojlim_s \CE^\bullet(\mathfrak g,
      D^{*,[r,s]})=\CE^\bullet(\mathfrak g, \Drigdeath(V^*)_K(\chi_{LT}))\] as
    $\mathfrak g$ is finite-dimensional.
    Analogously we also get a morphism \[ C^\bullet(\Gamma,\Drigdeath(V)_K)\ra
      \CE^\bullet(\mathfrak g, \Drigdeath(V)_K).\]
    By 
    \ref{prop:identify-duals-over-robba,prop:robba-dual-is-top-dual}, we can identify
    $\Drigdeath(V^*)_K(\chi_{LT})$ with
    $(\Drigdeath(V)_K)'_b=\Hom_{K,\mathrm{cts}}(\Drigdeath(V)_K,K)_b,$ where the strong dual is now
    taken over $K$. By \ref{prop:chevalley-eilenberg-duality} we get a
    $\Gamma$-equivariant $K$-linear morphism \[ \CE^\bullet(\mathfrak g,
      \Hom_{K,\mathrm{cts}}(\Drigdeath(V)_K,K)) \ra
      \Hom_{K,\mathrm{cts}}(\CE^\bullet(\mathfrak g, \Drigdeath(V)_K), K)[-1].\]

    Composing all these morphism, we get a functorial morphism of complexes \[
      C^\bullet(\Gamma, \Drigdeath(V^*)_K(\chi_{LT})) \ra
      \Hom_K(C^\bullet(\Gamma,\Drigdeath(V)_K), K)[-1],\] which we can extend to a
    double complex as follows:
    \begin{equation}
      \label[diagram]{diag:herr-duality}
      \tag{$\star$}
      \begin{tikzcd}[column sep=large]
        C^\bullet(\Gamma,\Drigdeath(V^*)_K(\chi_{LT})) \arrow{r}\arrow{d}{C(G,\varphi'-1)} & \Hom_K(C^\bullet(\Gamma,\Drigdeath(V)_K),K)[-1]\arrow{d}{\Hom_K(C(\Gamma,\varphi-1),K)}\\
        C^\bullet(\Gamma,\Drigdeath(V^*)_K(\chi_{LT})) \arrow{r} & \Hom_K(C^\bullet(\Gamma,\Drigdeath(V)),K)[-1]
      \end{tikzcd}
    \end{equation}
    Here, $\varphi$ denotes the intrinsic $\varphi$-operator on $\Ddeath(V)_K$
    and $\varphi'$ its vector space dual. Note that the dualised $\varphi$
    operator on $\Hom_{\Brigdeath{L}}(\Drigdeath(V),\Brigdeath{L})$ is the
    intrinsic $\psi$-operator on $\Drigdeath(V^*)$ and vice versa,
    cf.~\cite[remarks 4.7 and 5.6]{MR3629658}. The diagram can hence also be
    written as \[
      \begin{tikzcd}[column sep=large]
        C^\bullet(\Gamma,\Drigdeath(V^*)_K(\chi_{LT})) \arrow{r}\arrow{d}{C(G,\psi-1)} & \Hom_K(C^\bullet(\Gamma,\Drigdeath(V)_K),K)[-1]\arrow{d}{\Hom_K(C(\Gamma,\varphi-1),K)}\\
        C^\bullet(\Gamma,\Drigdeath(V^*)_K(\chi_{LT})) \arrow{r} & \Hom_K(C^\bullet(\Gamma,\Drigdeath(V)),K)[-1],
      \end{tikzcd}
    \]
    where $\psi$ is the intrinsic $\psi$-operator of
    $\Drigdeath(V^*)_K(\chi_{LT})$. 
    
    By \cite[theorem 11.6]{thomas:cohomology-of-topologised-monoids}, this induces the first
    morphism of complexes as required. The second morphism can be constructed
    completely analogously: Instead of using the intrinsic $\varphi$-operator of
    $\Ddeath(V)_K$ in \ref{diag:herr-duality} on the right hand side, start with
    its intrinsic $\psi$-operator. Then we get the vector space dual $\psi'$ on
    the left hand side, which is the intrinsic $\varphi$-operator of
    $\Ddeath(V^*)_K$.
  \end{proof}
\end{proposition}

\begin{remark}
  The comparison morphism \[ C^\bullet(\Gamma,\Drigdeath(V)_K)\ra \CE^\bullet(\mathfrak g,\Drigdeath(V)_K)\]
  does probably not induce an isomorphism on cohomology after taking
  $G$-invariant on the right hand side. We expect $\varprojlim^1$-terms to
  appear. Note however that for the first cohomology group, a Mittag-Leffler
  argument makes a comparison possible, cf.~\cite[proposition 2.1.1]{MR3690276}.

  Apart from this and under the assumptions of
  \ref{prop:duality-morphism-analytic-cohomology}, i.\,e., strict differentials
  in the Chevalley-Eilenberg complex and an open subgroup of $\Gamma$ operating
  trivially on the Lie algebra cohomology, we can follow the same argument to
  compare cohomology groups, as by
  \ref{prop:phi-gamma-modules-are-of-countable-type,prop:countable-type-exact-dual}
  taking duals is exact.
\end{remark}

\begin{remark}
  In degrees zero and one, $\varphi$ and $\psi$ yield the same cohomology
  groups, cf.~\cite[corollary 2.2.3]{MR3690276}.
\end{remark}

\printbibliography
\vfill
{\hfill\footnotesize The \textsc{Bib}\TeX-entries for this bibliography were
  mostly taken from MathSciNet.}

\end{document}